\documentclass[12pt]{article}
\usepackage[cp1251]{inputenc}
\usepackage{a4wide}
\usepackage{amsmath}
\usepackage{amssymb,bbm,verbatim}

\newcommand{\be}{\begin{equation}}
\newcommand{\ee}{\end{equation}}
\newcommand{\ba}{\begin{eqnarray}}
\newcommand{\ea}{\end{eqnarray}}
\newcommand{\ban}{\begin{eqnarray*}}
\newcommand{\ean}{\end{eqnarray*}}

\newtheorem{theo}{Theorem}

\newtheorem {prop} [theo] {Proposition}
\newtheorem {coro} [theo] {Corollary}

\newtheorem {rem} [theo] {Remark}

\newcommand{\supp}{\mathop{\rm supp}\nolimits}

\newcommand{\pfi}{\mathop{\varphi}\nolimits}

\newcommand{\N}{\mathop{\Bbb N}\nolimits}

\newcommand{\ZZ}{\mathbb Z}
\newcommand{\RR}{\mathbb{R}}
\newcommand{\NN}{\mathbb{N}}

\newcommand{\abcd}[2]{\hbox to\textwidth{#1\dotfill #2}}

\title{Wavelet frames on Vilenkin groups and their approximation properties
}
\author{
Yu.A. Farkov\footnote{Department of Applied Information Technology, Russian Presidential Academy of National Economy and Public Administration, Prospect Vernadskogo, 82, Moscow, 119571, Russia},
E.A. Lebedeva\footnote{Mathematics and Mechanics Faculty, Saint Petersburg State University,
Universitetsky prospekt, 28, Peterhof,  Saint Petersburg,
 198504, Russia; Saint Petersburg State Polytechnical University,
 Polytechnicheskay 29,  Saint Petersburg, 195251, Russia},
 M.A. Skopina\footnote{Faculty of Applied Mathematics and Control Processes, Saint Petersburg State University,
Universitetsky prospekt, 35, Peterhof, Saint Petersburg,   198504,
Russia}
}
\date{
farkov@list.ru, 
ealebedeva2004@gmail.com, 
skopina@ms1167.spb.edu}

\begin{document}
\maketitle

\setcounter{page}{1}


\begin{abstract}


 An explicit description of all Walsh polynomials generating tight wavelet frames is  given. 
An algorithm for finding the corresponding wavelet functions is suggested, and a general form 
for all wavelet frames  generated by an appropriate Walsh polynomial is described. Approximation  properties of tight wavelet frames 
are also studied.   In contrast to the real setting, it appeared that a wavelet tight frame
decomposition  has an arbitrary large approximation order whenever all wavelet functions are compactly supported.
\end{abstract}



\section{Introduction}
\label{s1}

In the forties of the last century Vilenkin started to study harmonic analysis on the groups which are infinite products of cyclic groups with the product topologies (see~\cite{1,2}). A well-known  special case is the locally compact 
Cantor group which is the countably weak direct product of the cyclic group of order 2. The foundations of harmonic analysis on the Vilenkin/Cantor groups (or the Walsh analysis) are presented in monographs~\cite{2},~\cite{3},~\cite{4}.
{\sloppy

}

The  wavelet theory for Vilenkin groups  is actively studied for the last two decades. Most of the results are related to the  Vilenkin groups  which are the products of cyclic groups of the same order $p$.
In the present paper we also discuss only such  groups.
 In 1996-1998 the concept  of   multiresolution analysis (MRA) for the Cantor  group was introduced by Lang  who also developed
a general method for the construction of MRA-based orthogonal wavelet bases (see~\cite{5} - \cite{7}). Some wavelets on the Cantor group were studied independently by Sendov~\cite{8} in 1997.
Later  deep investigations of orthogonal wavelets were added by one of the authors and his coauthors in~\cite{9} - \cite{13}.
An MRA theory and methods for the construction of orthogonal and biorthogonal wavelets for Vilenkin groups were developed by the same author jointly with Rodionov in~\cite{14}, \cite{15}.  Recently some results were added by Lukomskii in~\cite{16} 
and by one of the authors jointly with Krivoshein in~\cite{27}.

The wavelet theories for the  Vilenkin groups and for the real line are more or less similar. To construct a real compactly supported orthogonal wavelet basis one starts with an appropriate trigonometric polynomial (scaling mask) $m_0$. 
A scaling mask generating an orthogonal basis satisfies the condition  $|m_0(x)|^2+|m_0(x+\pi)|^2=1$ which is only necessary but not sufficient for the orthogonality. 
In the case of Vilenkin groups the  Walsh polynomials, i.e. finite linear combinations of the characters, play the same role as the trigonometric polynomials in the real setting. 
A necessary condition for a scaling mask to provide orthogonality is known. A complete
description of  Walsh polynomials satisfying this condition (scaling masks)  is given in~\cite{9}. Not all such scaling masks lead to orthogonal wavelet bases. Sufficient conditions are also known ~\cite{10, 13},
 but not for all masks  it is possible to check them easily without  computer calculations. However, similarly to the real setting, each such a scaling mask leads to a tight wavelet frame.
 Moreover, the necessary condition for the orthogonality can be relaxed in the case of the tight frame constructions.  Several examples of tight wavelet frames on the Cantor group are given in~\cite{21}.

The goals of this paper are the following. We give an explicit description of all Walsh polynomials
generating tight wavelet frames. We give an algorithmic method for finding the corresponding wavelet masks (matrix extension problem) and give an explicit description of all solutions to this matrix extension problem. 
 Note that there was an attempt to solve the latter problem in~\cite{17}. Unfortunately,
the attempt was unsuccessful, the suggested method does not lead to a required
matrix extension (see more detailed comments in Section~\ref{s3}). Note also that
 the matrix extension problem for the Vilenkin groups is not similar to the  same problem in the real setting, and our method is not an analog of the result obtained  by Lawton,
 Lee and Shen~\cite{18} for the trigonometric polynomials. Finally, we study approximation  properties of tight wavelet frame decompositions generated by  Walsh polynomials.
 It appeared that, in contrast to the real setting, any such decomposition has an arbitrary  approximation order.

\section{Notations and basic facts}
\label{s2}

As usual, by $\ZZ$, $\ZZ_+$, $\NN$, $\RR$, and $\Bbb C$ we denote the set of  integers, non-negative integers, positive integers, real numbers,
 and complex numbers, respectively.

Let an integer $p\ge2$ be fixed within the paper. The  Vilenkin group $G=G_p$ consists of the sequences
%
$$
 x=(x_j)=(\dots,0,0,x_k,x_{k+1},x_{k+2},\dots),
$$
where $x_j\in \{0,\dots,p-1\}$ for $j\in {\Bbb Z}$ and there exists at most a finite number of
negative $j$ such that  $x_j\ne0$.
The zero sequence is denoted by $\theta$.  If  $x\ne\theta$,  then there exists
a unique  $k=k(x)$ such that $x_k\ne0$ and  $x_j=0$ for  $j<k$.
The group operation on $G$  is denoted by~$\oplus$ and defined as the
coordinatewise addition modulo $p:$
$$
(z_j)=(x_j)\oplus (y_j) \, \Longleftrightarrow\, z_j = x_j + y_j \,({\rm
mod}\, p)\quad\mbox{for}\quad j\in {\Bbb Z}.
$$
Evidently, $\theta$ is the neutral element for $G$. 
We denote by $\ominus$  the inverse operation of $\oplus$. 
If $x\in G$, then $\ominus x$  denotes the inverse element of $x.$  

Define a mapping  \, $\lambda:\,G\to [0,+\infty)$ by letting
$$
    \lambda(x)=\sum_{j\in{\Bbb Z}}x_jp^{-j}, \qquad x=(x_j)\in G.
$$

The topology in $G$ is introduced via the complete system of neighborhoods of zero
$$
 U_l=\{(x_j)\in G\,|\ x_j=0 \quad   \mbox{for}\quad j\le l\}, \quad l\in {\Bbb Z}.
$$
Set $U=U_0$. Clearly, $U$ is a subgroup of $G$.
Introduce also the set
$$
H:=\{(x_j)\in G\,|\ x_j=0\ \mbox{for} \ j>0 \},
$$
 a  discrete subgroup of  $G$. It is clear that the image of $H$ under $\lambda$ is the set of non-negative
integers: $\lambda(H)={\Bbb Z}_+.$ For each $\alpha\in \mathbb{Z}_+$,  let  $h_{[\alpha]}$
denote the element of $H$ such that
$
\lambda(h_{[\alpha]})=\alpha;
$
in particular, $h_{[0]}=\theta$. A function $f$ defined on $G$ is called
{\em $H$-periodic} if $f(\omega\oplus h) = f(\omega)$ for all $h\in H$.

Define an automorphism $A$ on $G$  by letting  $(Ax)_j=x_{j+1}$.  Set
$$
       U_{n,\,s} := A^{-n}(h_{[s]})\oplus A^{-n}(U), \quad n\in{\Bbb Z}_+, \quad 0\le s\le p^n-1.
$$
It is clear that the sets $U_{n,\,s}$ are cosets of the subgroup $A^{-n}(U)$ in the group $U$.

There is a non-Archimedean metric on $G$ which generates the same topology. This metric is defined as follows.
The distance $\rho$ between $x\in G$ and $y\in G$ is given by $\rho(x,y)=\|x\ominus y\|$, where
 $\|\theta\|:=0$ and $\|x\|:= p^{-k(x)}$ for any $x\ne\theta$ (see~\cite[Section~9.1]{4}).

Since $G$ is a locally compact abelian group, there exists a Haar measure
$\mu$ on $G$ (see~\cite{HR}) which is  positive,
invariant under the shifts, i.e., $ \mu(x \oplus a)=\mu (x)$, and normalized by
$\mu(U)=1$.
Hence the functional spaces
$L_q(G)$ and $L_q(E)$, where $E$ is a measurable subset of $G$,
are defined. By $\mathbbm{1}_E$ we denote the characteristic function of the set $E$. Denote by $\langle\cdot\,,\cdot\rangle$  and  $||\cdot||_2$
the inner product and the norm in $L^2(G)$ respectively.

Given $\omega\in G$, the function
$$
\chi_\omega(x)=   \chi(x,\omega): = \exp\left(\frac{2\pi i}{p}\sum\limits_{j\in {\Bbb Z}} x_j\,\omega_{1-j}\right)
$$
is a group character of $G$. The Pontryagin dual group $G^*$ of $G$ is topologically isomorphic to $G$,
where the isomorphism is given by $\omega\to\chi_\omega$. In the sequel, we identify these groups and write $G$ instead of $G^*$.

It is easy to check the following properties of the characters

$$
\chi(x,\omega) = \chi(\omega,x),\ 
\chi(x_1\oplus x_2, \omega) = \chi(x_1,\omega)\chi(x_2,\omega),
$$
$$ 
\chi(A x,\omega) = \chi(x, A \omega),\ 
\overline{\chi(x,\omega)} = \chi(x,\ominus\omega).
$$

The Fourier transform of a function $f\in L^1(G)$ is defined by
$$
   \widehat f(\omega)=\int\nolimits_Gf(x)\overline{\chi (x,\omega)}d\mu(x), \quad \omega\in G.
$$
It is well-known that the Fourier operator
$$
                       {\cal F}\,:\ L^1(G)\cap L^2(G)\to L^2(G), \quad {\cal F}f=\widehat f,
$$
can be extended  in the standard way to the space $L^2(G)$ so that $  \langle\widehat f,\widehat g\,\rangle = \langle f,g \rangle$ for all  \,$f,g\in L^2(G)$\, ({\it Plancherel's relation}).

The {\it generalized Walsh functions} for the group $G$ are defined by
$$
    W_{\alpha}(x)=\chi(x,h_{[\alpha]}),  \qquad \alpha\in {\Bbb Z}_+,\, x\in G.
$$
These functions are $H$-periodic, continuous on $G$ and satisfy the orthogonality relations
$$
      \int_U W_{\alpha}(x) \overline{W_{\beta}(x)}d\mu(x)=\delta_{\alpha,\beta}, \qquad \alpha,\beta\in{\Bbb Z}_+,
$$
where  $\delta_{\alpha,\beta}$ is the Kronecker delta. Moreover,
the orthonormal system  $\{ W_{\alpha}\}$ is complete in $L^2(U)$
(see~\cite[Section 2.8]{3}). 
 If $0\le\alpha\le p^n-1$, then the function
 $ W_{\alpha}$ is constant on  each set $U_{n,\,s}$, $0\le s\le p^n-1.$
Finite linear combinations of generalized Walsh functions are called {\em Walsh polynomials}.

\medskip
{\it The discrete
Vilenkin-Chrestenson transform} maps a vector $\boldsymbol{b} = (b_0, \dots, b_{p^n-1})$ to the  vector  $\boldsymbol{a} = (a_0, \dots, a_{p^n-1})$ given by
\be
\label{1.1}
a_{\alpha}
=p^{-n}\sum_{s=0}^{p^n-1}b_sW_{\alpha}(A^{-n}h_{[s]}),
\quad 0\le\alpha\le p^n-1.
\ee
The inverse transform can be written as
\be
\label{1.2}
b_s=\sum_{\alpha=0}^{p^n-1}a_{\alpha}\,\overline{{W_{\alpha}(A^{-n}h_{[s]})}},
\qquad
0\le s\le p^n-1,
\ee
 see, e.g., ~\cite[Section 11.2]{3}  and~\cite{19} for the corresponding fast transforms.

\section{MRA-based tight wavelet frames}
\label{s3}

In this section we discuss  tight wavelet frames on the Vilenkin
groups constructed with using the MRA approach.
In the real setting, a general scheme for the construction of  MRA-based wavelet
frames was developed by Ron and Shen~\cite{20}. Following this scheme,
one starts with a refinable function (scaling function for the future MRA)
or its mask (scaling mask).
Then wavelet masks should be found according to the matrix extension
principle, and the wavelet masks provide the wavelet functions.
A similar approach is valid for  Vilenkin groups.

A  function $\pfi\in L^2(G)$ is called {\em refinable} if it satisfies
\be
\label{01}
         \varphi(x)=p\sum\limits_{\alpha=0}^\infty a_{\alpha}\varphi(Ax \ominus h_{[\alpha]}), \quad  x\in G.
\ee
To construct compactly supported wavelet frames we are interested in compactly supported refinable functions. In contrast to the real setting, the following statement holds true for the  Vilenkin groups.

\begin{prop}
If $\pfi$ is a compactly supported refinable function, then in~(\ref{01}) there exists only a finite number of nonzero coefficients $a_\alpha$.
\label{p1}
\end{prop}

{\bf Proof.} Let $\pfi$ satisfy~(\ref{01}) and $\supp\pfi\subset U_{-m}$,
$m\in\mathbb{Z}_+$. Set
 $$
 f(x)=p\sum_{\alpha\in\mathbb{Z}_+\atop \|\alpha\|\le p^{m+1}}a_{\alpha}
 \varphi(Ax \ominus h_{[\alpha]}),\ \ \
g(x)=p\sum_{\alpha\in\mathbb{Z}_+\atop \|\alpha\|> p^{m+1}}a_{\alpha}
\varphi(Ax \ominus h_{[\alpha]}).
$$
  If $x\in U_{-m}$ and $\|\alpha\|>p^{m+1}$, then
$\|Ax \ominus h_{[\alpha]}\|>p^{m+1}$, and hence $\varphi(Ax \ominus h_{[\alpha]})=0$. So, $g(x)=0$
whenever $x\in U_{-m}$.
If $x\not\in U_{-m}$ and $\|\alpha\|\le p^{m+1}$, then again
$\|Ax \ominus h_{[\alpha]}\|>p^{m+1}$, and hence $\varphi(Ax \ominus h_{[\alpha]})=0$.
 So, $f(x)=0$ whenever $x\not\in U_{-m}$. Since
also $\pfi(x)=0$ for all $x\not\in U_{-m}$, it follows from~(\ref{01})
that $\pfi=f$.
$\Diamond$

It is well known that in the real setting,  the mask of an orthogonal
compactly supported refinable function is a trigonometric polynomial.
An  analog of Proposition~\ref{p1} does not hold without additional
assumption of the orthogonality.

Let $\pfi$ be a compactly supported refinable function. Then, by Proposition~\ref{p1}, there exists a positive integer $n$ such that $\pfi$ satisfies a refinement equation
\be
         \varphi(x)=p\sum\limits_{\alpha=0}^{p^n-1}a_{\alpha}
				\varphi(Ax \ominus h_{[\alpha]}), \quad  x\in G. 
\label{1.6}
\ee
Taking the  Fourier transform of the both parts of (\ref{1.6}), we have
\be
\widehat{\pfi}(\omega)= m_0(A^{-1}\omega)\widehat{\pfi}(A^{-1}\omega),
\label{03}
\ee
where
\be
                m_0(\omega)=\sum\limits_{\alpha=0}^{p^n - 1}a_{\alpha}\overline{W_{\alpha}(\omega)}.   
                 \label{1.4}
\ee
The function $m_0$ is called {\em the mask} of $\pfi$ or  {\em  a scaling mask}.
This function is a Walsh polynomial of order $p^n-1$, $n\in\NN$. 
Indeed, using properties of a character of the group $G$, we get 
$\overline{W_{\alpha}(\omega)} = \overline{\chi(x,\,h_{[\alpha]})}=\chi(x,\,\ominus h_{[\alpha]})=
W_{\lambda(\ominus h_{[\alpha]})}(x)$. It follows from the definition of the operation $\oplus$ that $0\leq \alpha \leq p^n-1$ if and only if $0\leq \lambda(\ominus h_{[\alpha]}) \leq p^n-1.$  

Suppose there exist Walsh polynomials $m_1, \dots, m_r$
 $r\ge p-1$, ({\em wavelet masks}) such that  the matrix
\be
M(\omega)= \left(
\begin{array}{llll}
m_0(\omega)&m_1(\omega)&\ldots&m_r(\omega)
\\
m_0(\omega \oplus \delta_1)&m_1(\omega \oplus \delta_1)&\ldots&m_r(\omega \oplus \delta_1)
\\
\ \ \ \ \ldots&\ \ \ \ \ldots&\ldots&\ \ \ \ \ldots
\\
m_0(\omega \oplus \delta_{p-1})&m_1(\omega \oplus \delta_{p-1})&\ldots&m_r(\omega \oplus \delta_{p-1})
\end{array}
 \right), 
 \label{3.5}
\ee
where $\delta_l \in G$ and $\lambda(\delta_l)=l/p,$ $l\in\{0, \dots, p-1\}$,  satisfies 
\be
M(\omega)M^*(\omega)=I_p,
\label{100}
\ee
i.e., the rows of $M$ form an orthonormal system. By $M^*$ we denote the conjugate transpose of the matrix $M.$ 

The functions
 $\psi^{(1)},\dots, \psi^{(r)}$ defined by
$$
\widehat \psi^{(\nu)}(\omega)=m_{\nu}(A^{-1}\omega)\widehat\pfi(A^{-1}\omega), \ \ \,
\nu=1,\dots,r,
$$
are called {\em wavelet functions}. The corresponding {\em wavelet system}
\be
\psi^{(\nu)}_{j,k}(x) = p^{j/2}\psi^{(\nu)}(A^j x \ominus h_{[k]}), \quad  \nu=1,\dots, r, \ j\in \Bbb Z, \ k\in {\Bbb Z}_+. \label{1.9}
\ee
 is said to {\em be generated by $\pfi$}.

 A system $\{\psi^{(\nu)}_{j,k}\,|\  \nu=1,\dots,r, \ j\in \mathbb{Z}, \ k\in \mathbb{Z}_+\}$ is a {\it tight frame} in $L^2(G)$ if there exists a constant $C$
such that for any $f\in L^2(G)$ the following equality holds
$$
    \sum_{j\in \mathbb{Z}}\sum_{k\in \mathbb{Z}_+}\sum_{\nu=1}^r|\langle f\,,\psi^{(\nu)}_{j,k}\rangle|^2=C\|f\|^2.
$$
We say that a wavelet frame is compactly supported if each wavelet
function $\psi^{(\nu)}$ is compactly supported.

\begin{theo}[\cite{17}]
Let $\pfi$ be a compactly supported refinable function and $\widehat \pfi(\theta)\ne 0$. Then a wavelet system $\{\psi^{(\nu)}_{j,k}\}$ generated by $\pfi$ is a tight frame in $L^2(G)$ with $C=|\widehat \pfi(\theta)|^2$.
\label{t2.7}
\end{theo}

 An analog of this theorem for the real setting was proved by Petukhov~\cite{22}
(see also~\cite[\S~1.8]{23} or~\cite[\S~1.1]{24}).
The proof for the Vilenkin groups given in~\cite{17} is exactly the same. Moreover, following~\cite[\S~1.8]{23}, Theorem~\ref{t2.7} can be improved as follows.
 We can assume only the continuity of $\widehat\pfi$ at the point $\theta$ instead of the compactness of the support of $\pfi$.
 However, in the present paper we are interested only in compactly supported refinable functions.


Evidently, if the rows of matrix~(\ref{3.5}) form an orthonormal system, then
\be
\sum_{l=0}^{p-1}\left|m_0\left(\omega\oplus\delta_l\right)\right|^2  \le 1, \quad \ \omega \in G.
\label{1.7}
\ee
So, for a refinable function to generate a tight wavelet frame it is necessary
 its mask satisfy (\ref{1.7}). We claim that (\ref{1.7}) is also  sufficient,
i.e. if $m_0$ satisfies (\ref{1.7}),  then there exist wavelet masks $m_1, \dots, m_r$. To prove this and to give  an algorithmic method for the construction of wavelet masks (to solve  the matrix extension problem) we need several auxiliary statements.

\begin{prop}
For any $\alpha\in \mathbb{Z}_+$ and \,$k\in \{0,\dots, p-1\}$, we have
\be
\label{2.5}
                 W_{p\alpha + k}(\omega)=W_{\alpha}(A\omega)W_k(\omega), \quad \omega\in  G.
\ee
\label{p2}
\end{prop}

\textbf{Proof.}
Using the definition of the generalized Walsh function and properties of the character of $G$,
we immediately get
$$
     W_{p\alpha + k}(\omega)=\chi(\omega,h_{[p\alpha + k]})=\chi(\omega, h_{[p\alpha]} \oplus h_{[k]}) = \chi(\omega,h_{[p\alpha]})W_k(\omega),
$$
where
$$
      \chi(\omega,h_{[p\alpha]}) = \chi(\omega,Ah_{[\alpha]})= \chi(A\omega, h_{[\alpha]}) = W_{\alpha}(A\omega).\ \  \Diamond
$$

\medskip
\begin{prop}
For $\omega\in  G$, the matrix
\be
\label{2.6}
W(\omega)= \frac{1}{\sqrt p}\,\left(
\begin{array}{llll}
W_0(\omega)& W_1(\omega)&\ldots& W_{p-1}(\omega)
\\
W_0(\omega \oplus \delta_1)& W_1(\omega \oplus \delta_1)&\ldots& W_{p-1}(\omega \oplus \delta_1)
\\
\ \ \ \ \ldots&\ \ \ \ \ldots&\ldots&\ \ \ \ \ldots
\\
W_0(\omega \oplus \delta_{p-1})& W_1(\omega \oplus \delta_{p-1})&\ldots& W_{p-1}(\omega \oplus \delta_{p-1})
\end{array}
 \right),
 \ee
 where $\delta_l \in G$ and $\lambda(\delta_l)=l/p,$ $l\in\{0, \dots, p-1\}$,
is unitary.
\label{p3}
\end{prop}

\textbf{Proof.}
Suppose that $\omega\in  G$ and \ $k,s,s'\in \{0,\dots, p-1\}$. Then
$$
      W_k(\omega\oplus\delta_s)\overline{W_k(\omega\oplus\delta_{s'})}=W_k(\omega)W_k(\delta_s)\overline{W_k(\omega)}\,\overline{W_k(\delta_{s'})}
$$
$$
     =W_k(\delta_s\ominus \delta_{s'}) = \exp\left(\frac{2\pi i}{p}\,k(s - s')\right).
$$
Thus, we have
$$
     \sum_{k=0}^{p-1}W_k(\omega\oplus\delta_s)\overline{W_k(\omega\oplus\delta_{s'})}= p\delta_{s,s'}. \ \ \  \Diamond
$$

Fix $n\in\N$. Let
\be
\label{2.7}
                 m_0(\omega)=\sum\limits_{\alpha=0}^{p^n - 1}a_{\alpha}^{(0)}\overline{W_{\alpha}(\omega)}
\ee
be a Walsh polynomial satisfying~(\ref{1.7}). 
Then we fix $r\in\mathbb{N}$, 
$r\geq p-1$, \, and show how to find polynomials
\be
\label{2.8}
  m_{\nu}(\omega)=\sum\limits_{\alpha=0}^{p^n - 1}a_{\alpha}^{(\nu)}\overline{W_{\alpha}(\omega)}, \quad \nu=1,\dots, r,
\ee
such that the rows of matrix~(\ref{3.5}) form an orthonormal system.

\begin{prop} 
Let $m_0, \dots, m_r$ be polynomials given by~{\rm(\ref{2.7})} and~{\rm(\ref{2.8})}, and let
\be
\label{2.9}
  \mu_{\nu, k}(\omega):=\sqrt{p}\sum_{\alpha=0}^{p^{n-1} - 1}a_{p\alpha + k}^{(\nu)}\overline{W_{\alpha}(\omega)}, 
	\quad \nu=0,\dots, r, \ k=0,\dots, p-1.
\ee
 Then
\be
\label{2.15}
 m_{\nu}(\omega)
= \frac{1}{\sqrt p}\sum_{s=0}^{p-1}\mu_{\nu, s}(A\omega)\overline{W_{s}(\omega)}
\ee
and  for all  $s,s'\in \{0,\dots, p-1\}$ the following equivalence holds
\be
\label{2.10}
     \sum_{\nu=0}^{r}m_{\nu}(\omega\oplus\delta_s)\overline{m_{\nu}(\omega\oplus\delta_{s'})}= \delta_{s,s'} \ \Longleftrightarrow  \ \sum_{\nu=0}^{r}\mu_{\nu, s}(A\omega)\overline{\mu_{\nu, s'}(A\omega)}= \delta_{s,s'}.
\ee
Moreover, for all $\nu,\nu'\in \{0,\dots, r\}$, we have
\be
\label{2.11}
     \sum_{k=0}^{p-1}m_{\nu}(\omega\oplus\delta_k)\overline{m_{\nu'}(\omega\oplus\delta_k)}=\sum_{k=0}^{p-1}\mu_{\nu, k}(A\omega)\overline{\mu_{\nu', k}(A\omega)}.
\ee

\label{p4}
\end{prop}

{\bf Proof}. 
 Applying ~(\ref{2.5}) and using $H$-periodicity of the functions $W_{\alpha}$ for all $\nu\in \{0,\dots, r\}$ and $s\in \{0,\dots, p-1\}$, we obtain
$$
 m_{\nu}(\omega \oplus \delta_s)=\sum\limits_{\alpha=0}^{p^n - 1}a_{\alpha}^{(\nu)}\overline{W_{\alpha}(\omega \oplus\delta_s)}= \sum_{k=0}^{p-1}\sum_{\alpha=0}^{p^{n-1}-1}a_{p\alpha+k}^{(\nu)}\overline{W_{p\alpha+k}(\omega\oplus \delta_s)}
$$
\be
\label{2.12}
= \sum_{k=0}^{p-1}\sum_{\alpha=0}^{p^{n-1}-1}a_{p\alpha+k}^{(\nu)}\overline{W_{\alpha}(A\omega)} \, \overline{W_{k}(\omega \oplus\delta_s)}
= \frac{1}{\sqrt p}\sum_{k=0}^{p-1}\mu_{\nu, k}(A\omega)\overline{W_{k}(\omega \oplus\delta_s)}, 
\ee
which yields~(\ref{2.15}).  

Denote by $\Gamma(\omega)$ the matrix
$$
\Gamma(\omega)= \left(
\begin{array}{llll}
\mu_{0,0}(\omega)&\mu_{1,0}(\omega)&\ldots&\mu_{r,0}(\omega)
\\
\mu_{0,1}(\omega)&\mu_{1,1}(\omega)&\ldots&\mu_{r,1}(\omega)
\\
\ \ \ \ \ldots&\ \ \ \ \ldots&\ldots&\ \ \ \ \ldots
\\
\mu_{0,p-1}(\omega)& \mu_{1,p-1}(\omega)&\ldots&\mu_{r,p-1}(\omega)
\end{array}
 \right).
$$
By~(\ref{2.12}), $M(\omega)=\Gamma(A\omega)W^*(\omega)$,  which yields
$$
      M^*(\omega)M(\omega)= W(\omega)\Gamma^*(A\omega)\Gamma(A\omega)W^*(\omega).
$$
Therefore, we have
$$
      M^*(\omega)M(\omega)=I_{r+1}  \ \Longleftrightarrow  \ \Gamma^*(A\omega)\Gamma(A\omega)=I_{r+1},
$$
Thus~(\ref{2.10}) is proved. Finally, using Proposition \ref{p3} we get
$$
      M(\omega)M^*(\omega)= \Gamma(A\omega)(W^*(\omega)W(\omega))\Gamma^*(A\omega)= \Gamma(A\omega)\Gamma^*(A\omega),
$$
and~(\ref{2.11}) is proved.
$\Diamond$

The functions $\mu_{\nu, k}$ defined by~(\ref{2.9}) are called  {\em the polyphase components} of $m_\nu$.

\begin{coro}
Let  $m_0$ be defined  by  {\rm(\ref{2.7})} 
and satisfy~{\rm(\ref{1.7}),}  
and let $b^{(0,s)}_{n,l}:=\mu_{0,s}(A^{1-n}(h_{[l]}))$, where $\mu_{0, s}$ are the 
polyphase components of $m_0$. Then
\be
\label{2.16}
          \sum_{s=0}^{p-1}|b_{n,l}^{(0,s)}|^2\leq 1,\quad  l=0,\dots,p^n-1.
\ee
Moreover, the Walsh polynomials $m_\nu$, $\nu=1,\dots, r$, are wavelet functions corresponding to 
$m_0$ if and only if their polyphase components are given by 
\be
 \label{2.118}
 \mu_{\nu, s}(A\omega)=
 \sum_{l=0}^{p^n-1} b^{(\nu,s)}_{n,l}\mathbbm{1}_{U_{n,l}}(\omega),\quad \omega\in U,\quad  s=0,\dots, p-1,
\ee
where $b^{(\nu,s)}_{n,l}$ are complex numbers satisfying 
\be
\label{2.17}
          \sum_{\nu=0}^{r} b_{n,l}^{(\nu,s)}\overline{b_{n,l}^{(\nu,s')}} = \delta_{s,s'}, 
\quad  s,s'=0,\dots, p-1, \quad  l=0,\dots,p^n-1.
\ee
 \label{c1}
\end{coro}

{\bf Proof}. It follows from~(\ref{1.7}) and~(\ref{2.11}) that
$$
    \sum_{s=0}^{p-1}|\mu_{0, s}(A\omega)|^2 \leq 1, \quad \omega\in G,
$$
which yields~(\ref{2.16}).

By~(\ref{2.10}),   condition~(\ref{100})  is equivalent to
\be
\label{2.14}
\sum_{\nu=0}^{r}\mu_{\nu, s}(A\omega)\overline{\mu_{\nu, s'}(A\omega)}= \delta_{s,s'},
\quad \omega\in G, \quad  s,s'=0,\dots, p-1.
\ee
Since the polynomial  $\mu_{\nu, s}(A\cdot)$ is constant on the sets $U_{n,l}$, $l=0,\dots, p^n-1$, 
 we have~(\ref{2.118}). Using that $ \mathbbm{1}_{U_{n,l}}(\omega) \mathbbm{1}_{U_{n,l'}}(\omega)=\delta_{l,l'}$ for every
$\omega\in U$ and all $l, l'\in \{0,\dots p^n-1\}$, taking into account that $\mu_{\nu, s}$ is $H$-periodic and 
$A^{-n}(h_{[l]})\in {U_{n,l}}$,  we see that ~(\ref{2.14}) is equivalent to~(\ref{2.17}), 
as required. 
$\Diamond$

Next we must solve  system (\ref{2.17}) for each $l\in\{0,\dots, p^n-1\}$ with respect to
 the unknowns $ b_{n,l}^{(\nu,s)}$, $ \nu=1,\dots, r$, $s=0,\dots,p-1$, 
 provided (\ref{2.16}) is satisfied. It is easy to solve this problem using the following statement based on the Householder transform.

\begin{prop}
Let   $c_k\in  \mathbb{C}$,  $k=1,\dots, r$, $c_0\ne1$,  $\sum_{k=0}^r|c_k|^2=1$,  and let
\be 
c_{k0}=c_k,\quad
  c_{0j}=\overline{c}_{j0}\,\frac{1-c_{00}}{1-\overline{c}_{00}}, \quad c_{kj}= \delta_{kj} -  \frac{c_{k0}\overline{c}_{j0}}{1-\overline{c}_{00}},\quad j,k =1,\dots,r.
\label{HTr}
\ee
 Then the matrix \, $(c_{jk})_{j,k=0}^r$  \, is unitary.
\label{p5}
\end{prop}

\begin{rem}
We see that if $c_{0} = c_{00}$ is closed to $1,$ then calculations with the unitary matrix given by (\ref{HTr}) are numerically non-stable. However, to avoid this trouble it  suffices to renumber
properly  the initial vector $(c_k)_{k=0,\dots,r}.$
\end{rem}


Summarizing Propositions~\ref{p2}-\ref{p5}, we present an algorithm for the matrix extension. 
 on a given first column 
$(m_0(\omega ), \dots, m_0(\omega \oplus \delta_{p-1}))^T$, where $m_0$ is defined
 by  (\ref{2.7}) and satisfies~(\ref{1.7}).


\medskip
{\bf Algorithm  A.}


\begin{itemize}

\item
{\bf Step 1} Using (\ref{2.9}), find the polynomials $\mu_{0,s}$, $s=0,\dots p-1$.

\item
{\bf Step 2} Set $b^{*(0,s)}_{n,l}=\mu_{0,s}(A^{1-n}(h_{[l]}))$,  $s=0,\dots,p-1$, $l=0,\dots,p^n-1.$

\item
{\bf Step 3} In the case $r\ge p$,  for each $l$  choose arbitrary   $b^{*(0,p)}_{n,l}, \dots,b^{*(0,r)}_{n,l}\in \Bbb C$ such that
$\sum_{k=0}^{r}|b_{n,l}^{*(0,k)}|^2=1$ .

\item
{\bf Step 4} For each $l$, set  $c_k=b^{*(0,k)}_{n,l},$ $k=0,\dots, r$;
for each $s\in\{0,\dots,p-1\}$ and $\nu\in\{1,\dots,r\}$, compute 
  $c_{s,\nu}$ by~(\ref{HTr}) if
 $c_0\ne1$, and put $c_{s,\nu}=\delta_{s\nu}$ if $c_0=1$;    set $b^{*(\nu,s)}_{n,l}=c_{s,\nu}.$

\item
{\bf Step 5} Using (\ref{2.118}), find the polynomials $\mu_{\nu,s}$, $s=0,\dots p-1$, $\nu=1\dots,r$.

\item
{\bf Step 6} Using (\ref{2.15}), find $m_{\nu},$   $\nu=1\dots,r$.


As a result, we obtain a unitary matrix~(\ref{3.5}).

\end{itemize}

It is clear that the suggested method of matrix extension is not unique. First, there is a wide freedom for choosing 
 $b^{*(0,p)}_{n,l}, \dots,b^{*(0,r)}_{n,l}$  in Step 3. Second, Proposition~\ref{p5} provides only one unitary matrix on a given first column. 
 The following theorem gives a complete description of  all matrix extensions.

\begin{theo}
\label{t2}
Let  $m_0$ be defined  by  {\rm(\ref{2.7})} 
and satisfy~{\rm(\ref{1.7}),}   $b^{*(\nu,s)}_{n,l}$ be the numbers
determined  in Algorithm A, and let  $V_l:=(v^{(\nu,k)}_{n,l})_{\nu,k=0,\dots,r}$, $l=0,\dots,p^n-1,$ is a unitary matrix whose first  column is $(1,0,\dots,0)^{T}$. 
Then   polynomials $m_{\nu},$ $\nu=1,\dots,r$,  are wavelet masks corresponding to $m_0$ if an only if
$$
m_{\nu}(\omega)=\frac{1}{\sqrt{p}}\sum_{q=0}^{p^n-1}
\left(\sum_{s=0}^{p-1}c^{(\nu,s)}_{n,\lambda(h_{[q]}\ominus h_{[s]})}\right)
\overline{W_q(\omega)},
$$
where
$$
c^{(\nu,s)}_{n,t}=
\frac{1}{p^n}\sum_{l=0}^{p^n-1}
\sum_{k=0}^{r} b^{*(k,s)}_{n,l} v^{(\nu,k)}_{n,l}W_t(A^{1-n}(h_{[l]})) .
$$
\end{theo}
\textbf{Proof.}
It is well  known that  any  solution of  system~(\ref{2.17}) (for every fixed $l$) may be given as follows
$$
b^{(\nu,s)}_{n,l} = \sum_{k=0}^{r} b^{*(k,s)}_{n,l}{\tilde v}^{(\nu,k)}_{n,l},
$$
where  $\tilde V_l:=(\tilde v^{(\nu,k)}_{n,l})_{\nu,k=0,\dots,r}$ is a unitary matrix. Evidently,  to provide 
\be
b^{(0,s)}_{n,l}=b^{*(0,s)}_{n,l}, \quad s=0,\dots, p-1,
\label{200}
\ee
  the first  column of  $\tilde V$ should be equal to  $(1,0,\dots,0)^{T}$, 
that is
$$
b^{(\nu,s)}_{n,l} = \sum_{k=0}^{r} b^{*(k,s)}_{n,l}{ v}^{(\nu,k)}_{n,l}.
$$

Suppose Walsh  polynomials $m_{\nu},$ $\nu=1,\dots,r$,  given by~(\ref{2.8})  are wavelet masks corresponding to $m_0$. 
Let $ \mu_{\nu, s}$ denote the polyphase components of $m_\nu$, $\nu=0,\dots, p-1$.
Due to  Corollary~\ref{c1}, 
$$
 \mu_{\nu, s}(A\omega)= \sum_{l=0}^{p^n-1} b^{(\nu,s)}_{n,l} \mathbbm{1}_{U_{n,l}}(\omega),
$$
where  $ b^{(\nu,s)}_{n,l} $   satisfy~ (\ref{2.17})   for every $l$,  and the equality  $b^{(0,s)}_{n,l}=\mu_{0,s}(A^{1-n}(h_{[l]}))$ 
 is  equivalent to~(\ref{200}). 
Since the system $\{W_l\}_{l=0,\dots,p^n-1}$ is orthonormal on $U$ and each $W_l$ is constant on $U_{n,k},$
$k=0,\dots, p^n-1$,  it follows that
$$
 \mu_{\nu, s}(A\omega)= \sum_{l=0}^{p^n-1} c^{(\nu,s)}_{n,l}\overline{W_l(\omega)},
$$
where
 \ban
c^{(\nu,s)}_{n,t}=
\int_{U}
\left(\sum_{l=0}^{p^n-1} b^{(\nu,s)}_{n,l} \mathbbm{1}_{U_{n,l}}(\omega) \right)
W_t(\omega)\,d\omega=
\sum_{l=0}^{p^n-1} b^{(\nu,s)}_{n,l} \int_{U_{n,l}}
W_t(\omega)\,d\omega
\\
=\frac{1}{p^n}\sum_{l=0}^{p^n-1} b^{(\nu,s)}_{n,l} W_t(A^{1-n}(h_{[l]}))
=\frac{1}{p^n}\sum_{l=0}^{p^n-1}   \sum_{k=0}^{r} b^{*(k,s)}_{n,l}{ v}^{(\nu,k)}_{n,l}   W_t(A^{1-n}(h_{[l]})).
\ean
By (\ref{2.15}), we obtain
$$
m_{\nu}(\omega) =
\frac{1}{\sqrt{p}}  \sum_{s=0}^{p-1}\left(\sum_{l=0}^{p^n-1} c^{(\nu,s)}_{n,l}\overline{W_l(\omega)}\right)\overline{W_s(\omega)}
$$
Using the definition of the generalized Walsh function, of the function $h_{[\cdot]}$, and $\lambda$, we obtain
$$
W_l(\omega) W_s(\omega) = \chi(\omega,\,h_{[l]}) \chi(\omega,\,h_{[s]}) = \chi(\omega,\,h_{[l]}\oplus h_{[s]}) = W_{\lambda(h_{[l]}\oplus h_{[s]})}(\omega). 
$$  
We change index of summation. Denote 
$q = \lambda(h_{[l]}\oplus h_{[s]}),$
 then $h_{[q]}=h_{[l]}\oplus h_{[s]},$ so $h_{[l]}=h_{[q]}\ominus h_{[s]},$ thus $l=\lambda(h_{[q]}\ominus h_{[s]})$.
Finally,
$$
m_{\nu}(\omega) = \frac{1}{\sqrt{p}}\sum_{q=0}^{p^n-1}
\left(\sum_{s=0}^{p-1}c^{(\nu,s)}_{n,\lambda(h_{[q]}\ominus h_{[s]})}\right)
\overline{W_q(\omega)}.
$$

The converse follows from the same above  arguments. 
$\Diamond$

Another method for finding wavelet masks $m_1,\dots m_r$ was suggested
in~\cite[Theorem 3.5]{17}. Unfortunately, there are crucial mistakes there.
First,  the suggested form for the matrix
$\mathcal M_0(\xi) = \mathcal P(\xi)\sqrt{\Lambda(\xi)}\mathcal Q(\xi)$
 does not have required structure. Namely, the $k$-th row of the matrix $\mathcal M_0(\xi)$ should be $(m_1(\xi\oplus k/p), \dots, m_L(\xi\oplus k/p)).$ The author states that this is true for every unitary matrix $\mathcal Q(\xi)$. But
  if $\mathcal Q(\xi)$ is equal to the identity matrix,  some elements of the matrix $\mathcal P(\xi)\sqrt{\Lambda(\xi)}$ are equal to zero (for example, all elements of the first column except the first two elements), which is impossible. Second, even if the result were 
correct for a concrete  $\mathcal Q(\xi)$,  multiplying a matrix of the form (\ref{3.5})  by an arbitrary unitary matrix $\mathcal Q(\xi)$, one would lose the aforementioned structure.

\begin{theo}
A compactly supported refinable function generates compactly supported tight wavelet frames if and only if its mask $m_0$ is a Walsh polynomial satisfying $m_0(\theta)=1$
and~(\ref{1.7}).
\label{t1}
\end{theo}

{\bf Proof}. First we note that the assumption $\widehat \pfi(\theta)\ne 0$
in Theorem \ref{t2.7} is necessary for $\pfi$ to generate  tight wavelet frames.
Indeed, the proof of Theorem \ref{t2.7} is based on the identity
$$
\sum_{i=-\infty}^{j-1}\sum_{\nu=1}^r\sum_{k\in \mathbb{Z}_+}|\langle f\,,\psi^{(\nu)}_{i,k}\rangle|^2=\sum_{k\in \mathbb{Z}_+}|\langle f\,,\pfi_{j,k}\rangle|^2
$$
and  the relation
\be
\lim_{j\to+\infty}\sum_{k\in \mathbb{Z}_+}|\langle f,\pfi_{jk}\rangle|^2=|\widehat\pfi(\theta)|^2\|f\|^2,
\label{50}
\ee
which holds true for every $f\in L_2(G)$ Hence, if  $\{\psi^{(\nu)}_{j,k}\}$ is a tight frame, then  $\widehat \pfi(\theta)\ne 0$.   It follows from~(\ref{03}) that
$m_0(\theta)=1$. By Proposition \ref{p1} we conclude that $m_0$ is a Walsh polynomial. Property (\ref{1.7}) follows from the orthonormality of the rows of matrix (\ref{3.5}). Thus, necessity is proved.

If now $m_0$ is a Walsh polynomial and $m_0(\theta)=1,$ then 
$\widehat\pfi(\omega)=\prod_{j=1}^{\infty}m_0(A^{-j}\omega)$,  and hence, $\pfi(\theta)=1$. Existence of wavelet masks satisfying (\ref{3.5}) follows from  (\ref{1.7}) and Theorem~\ref{t2}.
$\Diamond$

\section{Generating scaling masks}
\label{s4}

In Section~\ref{s3} we suggested a method for the construction of wavelet masks 
$m_1,\dots, m_r$ (and hence of a tight wavelet frame) if an appropriate $m_0$ is given. 
Remind that $m_0$ should be a scaling mask, i.e. the mask of a compactly supported refinable function. In the present section, we describe all
such masks,  and hence describe all compactly supported  refinable functions generating tight wavelet frames.  Moreover, we give 
explicit formulas for these refinable functions.



The following statement is a minor improvement for an analog of  well-known Mallat's theorem
 (see, e.g.,  \cite[Lemma 4.1.3]{23}).

\begin{theo}
Suppose $m_0$ is a Walsh polynomial satisfying
 $m_0(\theta)=1$  and~(\ref{1.7}). Then the function
\be
           g(\omega):=\prod_{j=1}^{\infty} m_0(A^{-j}\omega), \quad \omega\in  G,
\label{2.2}
\ee
is in $L^2(G)$, and the function $\pfi$ defined by $\widehat\pfi=g$ is a compactly
supported refinable function whose mask is $m_0$.
\label{t4.1}
\end{theo}

{\bf Proof}.
Let $m_0$  be a Walsh polynomial of order $p^n-1$, $n\in\NN$. 
Since \, $m_0(\theta)=1$, we have \, $m_0(\omega)=1$ \, for all  \,$\omega\in U_{n,0}$. Hence for each $\omega\in G$, 
 at most a finite number of  factors of the product in~(\ref{2.2}) 
 are 
 not equal to $1$.  It follows that the product converges uniformly 
on any compact set  $K\subset G$ (because  $A^{-j}K\subset U_{n,0}$   whenever $j$ is large enough).

Set $g_0:=\mathbbm{1}_U$, 
$$
    g_l(\omega) := \prod_{j=1}^l  m_0(A^{-j}\omega) \mathbbm{1}_U(A^{-j}\omega), 
\quad \omega\in G, \quad l\in \mathbb{N}.
$$
Clearly we have
\be
             \lim_{l\to\infty}g_l(\omega) =g(\omega),          \quad \omega\in G,  
\label{4.3}
\ee
and
\be
             g_l(\omega) =  m_0(A^{-1}\omega)g_{l-1}(A^{-1}\omega), \quad  l\in \mathbb{N}. 
\label{4.4}
\ee
Consider now the $H$-periodic functions
$$
                 \Phi_l(\omega):=\sum_{h\in H}|g_l(\omega \oplus h)|^2,   \quad   l\in \mathbb{Z}_+.
$$
Note that if $k\in \{0,\dots,p-1\}$, $s\in \mathbb{Z}_+$, then $\delta_k=A^{-1}h_{[k]}$, $h_{[s]}=A^{-1}h_{[sp]}$,
$h_{[sp+k]}=h_{[sp]}\oplus h_{[k]}$. Therefore
$$
           W_{\alpha}(A^{-1}(\omega\oplus h_{[sp+k]}))= W_{\alpha}(A^{-1}\omega\oplus \delta_k)
$$
that,  by~(\ref{4.4}),  yields
$$
    \Phi_l(\omega)= \sum_{k=0}^{p-1}\sum_{s\in \mathbb{Z}_+}|g_l(\omega \oplus h_{[sp+k]})|^2
$$
$$
        = \left(|m_0(A^{-1}\omega)|^2 + |m_0(A^{-1}\omega\oplus\delta_1)|^2 + \dots +  |m_0(A^{-1}\omega\oplus\delta_{p-1})|^2 \right)\cdot \Phi_{l-1}(\omega).
$$
It follows from~(\ref{1.7}) that
$$
           |\Phi_l(\omega)| \leq |\Phi_{l-1}(\omega)|.
$$
Since $|\Phi_0(\omega)| \leq 1$,  we have $|\Phi_l(\omega)| \leq 1$ for all $l\in \mathbb{N}$, $\omega\in G$. Consequently,
$$
\int_G|g_l(\omega)|^2\,d\mu(\omega) = \sum_{h\in H}\int_{U\oplus h}|g_l(\omega)|^2\,d\mu(\omega)
$$
$$
  = \int_U\sum_{h\in H}|g_l(\omega\oplus h)|^2\,d\mu(\omega) =\int_U\Phi_l(\omega)\,d\mu(\omega)\leq 1, \quad l\in \mathbb{N}.
$$
According to~(\ref{4.3}) and Fatou's lemma, we then obtain
$$
    \int_G|g(\omega)|^2\,d\mu(\omega) \leq \sup_{l\in\ZZ_+} \int_G|g_l(\omega)|^2\,d\mu(\omega) \leq 1.
$$
Thus  $g\in L^2(G)$, and  hence, by  Plancherel's theorem,  $\pfi\in L^2(G)$.
Evidently, $\pfi$ satisfies~(\ref{03}), so  $\pfi$ is a refinable function   with the mask $m_0$. To complete the proof  it remains to show that $\pfi$ is compactly supported.
But this easily  follows from~\cite[Proposition~2]{13}. $\Diamond$

\begin{rem} Evidently, any function $\pfi$ defined by $\widehat\pfi=cg$, where $g$ is   
from Theorem~\ref{t4.1}, $c\ne0$, 
is  a compactly supported refinable function with the same mask $m_0$, and $\widehat\pfi(\theta)=c$. It follows from
 Theorems~\ref{t1} and ~\ref{t2.7}  that such $\pfi$ generates  tight wavelet frames whose frame boundary equals $|c|^2$.
\end{rem}

We next describe all Walsh polynomials satisfying assumptions of Theorem~\ref{t4.1}.

\begin{prop}
\label{p10}
Let  $m_0$ be  a Walsh polynomial of order $p^n-1$, $n\in\NN$, given by~(\ref{1.4}).
 Conditions  $m_0(\theta)=1$  and~(\ref{1.7}) hold true for the polynomial $m_0$  
  if and only if its coefficients  $a_\alpha$ are defined  by the
 Vilenkin-Chrestenson transform~(\ref {1.1}), where the numbers $b_s$ satisfy
\be
  b_0=1, \quad  |b_l|^2 +|b_{l + p^{n-1}}|^2 + \dots +|b_{l + (p - 1)p^{n-1}}|^2 \leq 1,  \quad 0\le  l\le p^{n-1}-1.
\label{4.5}
\ee
\end{prop}

{\bf Proof}. Relations~(\ref {1.1}) and~(\ref{1.2}) are equivalent.     It follows  from~(\ref{1.2}) that
 $b_s=m_0(A^{-n}h_{[s]})$,  $s=0,\dots, p^n-1$. Since  $  m_0(A^{-n}h_{[s]})=m_0(\omega)$ 
 for all $\omega\in U_{n,\,s}$ and $m_0$ is $H$-periodic,  
we see that  both conditions $m_0(\theta)=1$  and~(\ref{1.7}) are valid if and only if~(\ref{4.5}) is satisfied. $\Diamond$

Analyzing the latter statement, we can give an algorithmic method  for the construction of   Walsh polynomials satisfying
 assumptions of Theorem~\ref{t4.1}.

\medskip 
{\bf Algorithm  B.}

\begin{itemize}

\item
{\bf Step 1.}  Choose $n\in\NN$ and arbitrary  numbers $b_0,\dots, b_{p^n-1}\in\Bbb C$ 
satisfying~(\ref{4.5}).
{\sloppy

}

\item
{\bf Step 2.} 
Compute  $a_0,\dots, a_{p^n-1}$ by~(\ref {1.1}).

\item
{\bf Step 3.} Define $ m_0 $ by~(\ref {1.4}).

\end{itemize}

\section{Approximation order}
\label{s5}

In this section we study approximation properties of compactly supported wavelet frames. An MRA-based method for the construction of such frames was described in the previous sections, but here we consider arbitrary compactly supported wavelet frames, not necessary MRA based as above.

For  $m\in {\Bbb Z}_+$, we set
\be
W_2^m=\left\{\,f\in L_2(G)\,:\ \int\limits_G\|\omega\|^{2k}|\widehat{f}(\omega)|^2\,d\mu(\omega)<\infty \quad \mbox{for}\quad k=0,1,\dots, m \,\right\}.
\ee
It is not difficult to see that $W_2^m$ is a Banach space with the norm
$$
     \|f\|_{2,m}:=\sum_{k=0}^m\left(\int\limits_G\|\omega\|^{2k}|\widehat{f}(\omega)|^2\,d\mu(\omega)\right)^{1/2}.
$$
So, $W_2^m$ is an analog of the Sobolev space in the real setting.

\begin{theo}
Let  $f\in W_2^m$, $j\in\Bbb Z$, and let $\{\psi^{(\nu)}_{i,k}|\,   \nu=1,\dots, r, \ i\in \Bbb Z, \ k\in {\Bbb Z}_+,\}$ be a compactly supported tight wavelet frame. Then
\be
\Big\|f - \sum_{i=-\infty}^{j}\sum_{\nu=1}^r\sum_{k\in \mathbb{Z}_+}\langle f\,,\psi^{(\nu)}_{i,k}\rangle \psi^{(\nu)}_{i,k}\Big\|_2\le C p^{-jm}\|f\|_{2,m},
\label{02}
\ee
where $C$ does not depend of $f$ and $j$.
\label{t0}
\end{theo}

{\bf Proof}. Due to the main property of tight frames in
Hilbert spaces (see, e.g.,~\cite[\S 1.8]{23}),
every $f\in L_2(G)$ can be decomposed as
$$
f=\sum_{\nu=1}^r\sum_{i=-\infty}^{+\infty}\sum_{k\in \mathbb{Z}_+}\langle f\,,\psi^{(\nu)}_{i,k}\rangle \psi^{(\nu)}_{i,k},
$$
where the series converges unconditionally.
Hence
\ban
\Big\|f - \sum_{\nu=1}^r\sum_{i=-\infty}^{j}\sum_{k\in \mathbb{Z}_+}\langle f\,,\psi^{(\nu)}_{i,k}\rangle\psi^{(\nu)}_{i,k}\Big\|_2 = \Big\|\sum_{\nu=1}^r\sum_{i=j+1}^{\infty}\sum_{k\in \mathbb{Z}_+}\langle f\,,\psi^{(\nu)}_{i,k}\rangle \psi^{(\nu)}_{i,k} \Big\|_2
\\
\leq \sum_{\nu=1}^r\sum_{i=j+1}^{\infty} \Big\|\sum_{k\in \mathbb{Z}_+}\langle f\,,\psi^{(\nu)}_{i,k}\rangle\psi^{(\nu)}_{i,k}\Big\|_2
\ean
Since $\{\psi^{(\nu)}_{i,k}\}$ is a frame, the functions
$\psi(\omega\oplus h)$, $h\in H$, form a Bessel system.
It follows that
\be
\sum_{h\in H}|\widehat{\psi}(\omega\oplus h)|^2\le C_1\ \   \mbox{a.e.}
\label{04}
\ee
This fact can be easily checked in the same way as it is done for the real case
(see, e.g.,~\cite[Theorem 1.1.7]{23}).

For $i>j$ and $\psi=\psi^{(\nu)}$, using  the Plancherel theorem and~(\ref{04}),
we have
\ba
\Big\|\sum_{k}\langle f\,,\psi_{i,k}\rangle\psi_{i,k}\Big\|_2^2=\Big\|\sum_{k}\langle f\,,\psi_{i,k}\rangle\widehat{\psi}_{i,k}\Big\|_2^2
\nonumber
\\
=p^{-j}\int\limits_G\Big|\sum_{k}\langle f\,,\psi_{i,k}\rangle
W_k(A^{-i}\omega)\widehat{\psi}(A^{-i}\omega)\Big|^2\,d\mu(\omega)
\nonumber
\\
=\int\limits_G\Big|\sum_{k}\langle f\,,\psi_{i,k}\rangle W_k(\omega)\Big|^2 |\widehat{\psi}(\omega)|^2\,d\mu(\omega)
\nonumber
\\
=\int\limits_U\sum_{h\in H}|\widehat{\psi}(\omega\oplus h)|^2\,\Big|\sum_{k}\langle f\,,\psi_{i,k}\rangle W_k(\omega)\Big|^2\,d\mu(\omega)
\nonumber
\\
\leq C_1 \int\limits_U\Big|\sum_{k}\langle f\,,\psi_{i,k}\rangle W_k(\omega)\Big|^2\,d\mu(\omega)=
C_1 \sum_{k}|\langle f\,,\psi_{i,k}\rangle |^2.
\label{05}
\ea

It is easy to see that the function
$$
G_j(\omega):=p^i\sum_{h\in H}\widehat{f}(A^i(\omega\oplus h))\overline{\widehat{\psi}(\omega\oplus h)}
$$
is $H$-periodic, $G_j\in L_2(U)$ and
$$
            \langle f\,,\psi_{i,k}\rangle = \int\limits_U G_j(\omega)W_k(\omega)\,d\mu(\omega).
$$
It follows that
\be
\left(\sum_{k}|\langle f\,,\psi_{i,k}\rangle |^2\right)^{1/2}=
   \left(\int\limits_Up^i|\sum_{h\in H}\widehat{f}(A^i(\omega\oplus h))\widehat{\psi}(\omega\oplus h)|^2\,d\mu(\omega)\right)^{1/2}\le I_1+I_2,
   \label{06}
\ee
where
$$
I_1^2=\int\limits_Up^i|\widehat{f}(A^i(\omega))\widehat{\psi}(\omega)|^2\,d\mu(\omega), \quad
I_2^2=\int\limits_Up^i|\sum_{h\neq \theta}\widehat{f}(A^i(\omega\oplus h))\widehat{\psi}(\omega\oplus h)|^2\,d\mu(\omega).
$$
If $\omega\in U$, $h\in H$ and $h\ne\theta$, then $\|\omega\oplus h \|=\|h\|\ge1$.
Hence, by~(\ref{04}),
\ba
I_2^2\leq p^i\int\limits_U\sum_{h\neq \theta}\|\omega\oplus h \|^{2m}|\widehat{f}(A^i(\omega\oplus h))|^2\cdot
\sum_{h\neq \theta}\frac{|\widehat{\psi}(\omega\oplus h)|^2}{\|\omega\oplus h \|^{2m}}\,d\mu(\omega)
\nonumber
\\
\leq C_1p^i\int\limits_G\|\omega\|^{2m}|\widehat{f}(A^i\omega)|^2\,d\mu(\omega)\leq C_1p^{-2im}\|f\|^2_{2,m}.
\label{07}
\ea

Since $\psi$  is a compactly supported function,
 there exists a positive integer $n$ such that $\widehat\psi$  is constant 
on each set $U_{n,s}$ 
 (see, e.g.,~\cite[Proposition~2]{13}). On the other hand,  $\widehat\psi$ is continuous  and        the series
$ \sum_{i=0}^{\infty}\sum_{k\in \mathbb{Z}_+}|\langle f\,,\psi_{i,k}\rangle|^2$  is convergent.
It  follows from~(\ref{50}) with $\pfi=\psi$ that $\widehat\psi(\theta)=0$. 
Thus  $\widehat\psi(\omega) = 0$ for all $\omega\in U_{n,\,0}$. Using this and
taking into account that $\|\omega\|\ge p^{-n}$  whenever $\omega\not\in U_{n,\,0}$, we have
\ba
I_1^2 = p^i\int\limits_{U\setminus U_{n,\,0}}|\widehat{f}(A^i\omega)\widehat{\psi}(\omega)|^2\,d\mu(\omega)
\leq \frac{\|\widehat\psi\|_\infty\,p^i}{p^{-2mn}}\int\limits_{U\setminus U_{n,\,0}}\|\omega\|^{2m}
|\widehat{f}(A^i\omega)|^2\,d\mu(\omega)
\nonumber
\\
\leq \frac{\|\widehat\psi\|_\infty\,p^i}{p^{-2mn}}
\int\limits_G\|\omega\|^{2m}|\widehat{f}(A^i\omega)|^2\,
d\mu(\omega)\leq C_2 p^{-2im} \|f\|^2_{2,m}.
\label{08}
\ea
It follows from~(\ref{06})-(\ref{08}) combined with (\ref{05}) that
$$
\Big\|\sum_{k}\langle f\,,\psi^{(\nu)}_{i,k}\rangle\psi^{(\nu)}_{i,k}\Big\|_2=
\Big\|\sum_{k}\langle f\,,\psi_{i,k}\rangle\psi_{i,k}\Big\|_2\le C_3p^{-im} \|f\|_{2,m},
$$
and hence
$$
\sum_{\nu=1}^r\sum_{i=j+1}^{\infty} \Big\|\sum_{k\in \mathbb{Z}_+}\langle f\,,\psi^{(\nu)}_{i,k}\rangle\psi^{(\nu)}_{i,k}\Big\|_2\le C p^{-jm}\|f\|_{2,m},
$$
which was to be proved. $\Diamond$

Note that an analog of Theorem~\ref{t0} does not exist in the real setting. It is well known that some special property of  wavelet functions $\psi^{(\nu)}$  should be satisfied for the corresponding frame decomposition to provide approximation order $m$.
 We see that the situation is different in the case of Vilenkin groups. Wavelet tight frame
decomposition  has an arbitrary large approximation order whenever all wavelet functions are compactly supported.

\section*{Acknowledgments}
Yu.A.F.  is supported by the Ministry of Education and Science of the Russian Federation, grant
\#5.5805.2011.
E.A.L. is supported by the RFBR, grant \#12-01-00216, 
 by Saint Petersburg State University, grant  \#9.38.62.2012,
and by grant of President RF \#MK-1847.2012.1.
M.A.S. is supported by the RFBR-grant \#12-01-00216, and by Saint Petersburg State University, grant \#9.38.62.2012.

\end{document}